\documentclass[a4paper,12pt,leqno]{article}
\usepackage{amssymb,amsmath}
\usepackage[applemac]{inputenc}
\usepackage{xypic}
\xyoption{all}
\usepackage{graphicx}
\usepackage{hyperref}

\newtheorem{defn}{Definition}[section]

\newtheorem{corollary}[defn]{Corollary}
\newtheorem{rem}[defn]{Remark}
\newtheorem{exm}[defn]{Example}
\newtheorem{lemma}[defn]{Lemma}
\newtheorem{theorem}[defn]{Theorem}
\newtheorem{notat}[defn]{Notation}
\newtheorem{newpar}[defn]{}

\newtheorem{xdefn}{Definition.}
\newtheorem{xproposition}{Proposition.}
\newtheorem{xcorollary}{Corollary.}
\newtheorem{xrem}{Remark.}
\newtheorem{xexm}{Example.}
\newtheorem{xlemma}{Lemma.}
\newtheorem{xtheorem}{Theorem.}
\newtheorem{xnotat}{Notation.}
\newtheorem{xnewpar}{\it}
\newtheorem{xproof}{{\it Proof. }}
\newtheorem{xproofof}{{\it Proof}}

\newenvironment{definition}{\begin{defn}\em}{\end{defn}}
\newenvironment{remark}{\begin{rem}\em}{\end{rem}}

\newenvironment{proof}{\begin{xproof}\em}{\end{xproof}}

\newenvironment{newparagraph*}[1]{\begin{xnewpar}\hspace*{-1.5mm}{#1}. \rm}{\end{xnewpar}}

\newenvironment{definition*}{\begin{xdefn}\em}{\end{xdefn}}
\newenvironment{remark*}{\begin{xrem}\em}{\end{xrem}}
\newenvironment{example*}{\begin{xexm}\em}{\end{xexm}}
\newenvironment{notation*}{\begin{xnotat}\em}{\end{xnotat}}
\newenvironment{proposition*}{\begin{xproposition}}{\end{xproposition}}
\newenvironment{corollary*}{\begin{xcorollary}}{\end{xcorollary}}
\newenvironment{lemma*}{\begin{xlemma}}{\end{xlemma}}
\newenvironment{theorem*}{\begin{xtheorem}}{\end{xtheorem}}

\newenvironment{eq}{\setcounter{equation}{\arabic{defn}}\begin{equation}}{\end{equation}\setcounter{defn}{\arabic{equation}}}

\def\qed{\hspace{0.3cm}{\rule{1ex}{2ex}}}
\newcommand\V{\bigvee}

\newcommand\ie{i.e.}
\newcommand\eg{e.g.}

\newcommand\st{\mid}

\newcommand\topology{\operatorname{\Omega}}
\newcommand\groupoid{\operatorname{\mathcal{G}}}

\newcommand\downsegment{{\downarrow}}

\newcommand\ident{\mathrm{id}}
\newcommand\ipi{\mathcal I}

\newcommand\lcc{\operatorname{{\mathcal L}^{\vee}}}

\newcommand\germ{\operatorname{germ}}
\newcommand\cod{\operatorname{cod}}

\newcommand\dist{\operatorname{dist}}

\newcommand\Germs{\mathrm{Germs}}
\newcommand\id{\mathrm{id}}

\newcommand\dom{\operatorname{dom}}

\newcommand\hide[1]{}

\begin{document}

\title{\'Etale groupoids as germ groupoids and their base extensions\thanks{Research supported in part by Funda\c{c}\~{a}o para a Ci\^{e}ncia e a Tecnologia through FEDER and project PPCDT/MAT/55958/2004.}
}
\author{\sc Dmitry Matsnev and Pedro Resende
}
\date{~\vspace*{-1cm}}

\maketitle

\begin{abstract}
We introduce the notion of wide representation of an inverse semigroup and prove that with a suitably defined topology there is a space of germs of such a representation which has the structure of an \'etale groupoid. This gives an elegant description of Paterson's universal groupoid and of the translation groupoid of Skandalis, Tu, and Yu. In addition we characterize the inverse semigroups that arise from groupoids, leading to a precise bijection between the class of \'etale groupoids and the class of complete and infinitely distributive inverse monoids equipped with suitable representations, and we explain the sense in which quantales and localic groupoids carry a generalization of this correspondence.
\end{abstract}


\section{Introduction}

The local bisections of a topological \'etale groupoid $G$ can be identified with the open subsets of the space of arrows on which the domain and the range maps are injective. The set of all the local bisections of $G$ forms an inverse semigroup, denoted by $\mathcal I(G)$. This inverse semigroup is equipped with a natural representation on the space $G_0$ of units of $G$, that is a homomorphism to the inverse semigroup of partial homeomorphisms on the space of units. This representation is \emph{full} in the sense that the idempotents of $\mathcal I(G)$ bijectively correspond to the open sets of $G_0$.

On the other hand, for any inverse semigroup $S$ equipped with a full representation on a topological space $X$, its space of \emph{germs} with a sheaf-like topology can be given a structure of \'etale groupoid with $X$ as unit space. The first objective of this paper is to extend the aforementioned germ groupoid construction to the following more general case, where an inverse semigroup $S$ represented on a topological space $X$ will be called \emph{wide over $X$} if the images of the idempotents in $S$ under the representation cover $X$:

\begin{theorem*}
For any wide inverse semigroup over a topological space $X$, its space of germs can be given the structure of an \'etale groupoid with object space $X$.
\end{theorem*}

An immediate application of this theorem is the following extension procedure for \'etale groupoids: (1) start with an \'etale groupoid $G$ plus the full representation of $\ipi(G)$ on $G_0$; (2) define an extension $X$ of the space $G_0$ in such a way that $\ipi(G)$ becomes a wide inverse semigroup over $X$; (3) apply the above theorem in order to obtain an \'etale groupoid which is the required ``base extension'' of the original one.

The second objective of this paper is to describe how two classical constructions in the theory of \'etale groupoids can be achieved using the germ groupoid approach.

The first example is the \emph{universal groupoid} of Paterson, a certain \'etale groupoid associated to a countable inverse semigroup, with the properties that the $C^*$-algebras for the semigroup and its groupoid are the same, both the full and the reduced ones. The original construction~\cite{Paterson} relies on so-called localisation techniques due to \cite{Kumjian}; we show that it can be described more succinctly in terms of the above theorem, as the germ groupoid associated to a wide representation.

Concerning the second example, suppose we start with a discrete groupoid $G$ (this is the same as an \'etale groupoid whose unit space is discrete). Then the representation of the underlying inverse semigroup can be extended to a representation by partial homeomorphisms on the Stone--\v Cech compactification of the unit space of $G$. Taking the germs of this extension, we obtain an \'etale groupoid $\beta_0G$, which extends $G$ in such a way that its unit space is the Stone--\v Cech compactification of the unit space of $G$. The arrows of $\beta_0G$ are simply the arrows of $G$ and the other necessary ones required by the structure laws (hence $\beta_0G$ can be viewed as the most economical extension of $G$ with the required space of units).

This result is closely related to work of Skandalis, Tu, and Yu, who in~\cite{SkandalisTuYu} established the connection between \'etale groupoids and coarse metric spaces (see \cite{Roe}). They defined the \emph{translation groupoid} of a boundedly discrete coarse space as the maximal possible principal \'etale extension of the pair groupoid on that space such that the unit space of the extension is the Stone--\v{C}ech compactification of the space we started with. The coarse Baum--Connes conjecture for the original space is then equivalent to the Baum--Connes conjecture for its translation groupoid. As a consequence of that, for example, the Hilbert space embeddability for a discrete space implies the coarse Baum--Connes conjecture for such a space, and, via the Higson descent technique, one gets an elegant proof of the Novikov conjecture for discrete groups equipped with the Hilbert space embeddings.

If we take the inverse semigroup generated by controlled partial bijections of $X$ (see~\cite{Roe}) with its tautological representation on $X$ and extend it via the Stone--\v Cech technique described above, the resulting inverse semigroup represented on $\beta X$ is often wide, for instance when the coarse space $X$ is unital (\eg, this happens when the coarse structure comes from the equivalence class of metrics, which is particularly relevant for the Novikov conjecture). In such a case the translation groupoid can be obtained by our germ groupoid construction.

Finally, as a third objective we obtain a precise characterization of the inverse semigroups of the form $\ipi(G)$ and a bijection (up to isomorphisms) between the class of \'etale groupoids with unit space $X$ and the class of complete inverse semigroups equipped with full representations on $X$. The particular case where $X$ is a sober space leads, via \cite{Resende}, to a connection with localic groupoids (internal groupoids in the category of locales) and quantales, whose significance in the context of the present paper we also describe.

The paper is structured as follows. In Section~\ref{section:GermGroupoids} we recall the notions of inverse semigroup, inverse semigroup representation, \'etale groupoid, and germs. We continue with the structure theorems connecting these notions and then introduce wide representations. The main result of this section is the theorem stated above.
Section~\ref{section:UniversalGroupoid} is devoted to the discussion of the universal groupoid of Paterson.
In Section~\ref{section:TranslationGroupoid} we present the Stone--\v Cech extension technique for discrete groupoids. We review the translation groupoid of Skandalis, Tu, and Yu, and discuss how it arises as a germ groupoid.
Finally, in Section~\ref{section:ACPs} we close the circle by characterizing the inverse semigroups that arise from \'etale groupoids (as described in Section~\ref{section:GermGroupoids}). This requires us to address the order-theoretic properties of inverse semigroups and, in particular, the notion of complete inverse semigroup. We conclude with the remarks on localic groupoids and quantales that establish the relation to \cite{Resende}.

\section{Germ groupoids of inverse semigroups}\label{section:GermGroupoids}

\subsection*{Inverse semigroups}

\begin{definition}\label{def:inversesemigroup}
Let $S$ be a semigroup. An element $t\in S$ is said to be an \emph{inverse} of $s\in S$ if
\[sts = s \qquad\mbox{ and }\qquad tst=t.\]
An \emph{inverse semigroup} $S$ is a semigroup such that for all $s\in S$ there is a unique inverse of $s$, which we shall denote by $s^*$.
Equivalently, an inverse semigroup is a semigroup for which each element has an inverse and for which any two idempotents commute (see \cite{Lawson}). The set of idempotents is denoted by $E(S)$.
An \emph{inverse monoid} is an inverse semigroup that has a multiplicative unit, which we shall denote by $e$.
\end{definition}

A few simple but useful facts about inverse semigroups are:
\begin{enumerate}
\item In any inverse semigroup the inverse operation defines an involution.
\item The set of idempotents forms a semilattice with meet given by multiplication: $f\land g=fg$; this semilattice has a greatest element if and only if $S$ is a monoid.
\item The idempotents are precisely the elements of the form $ss^*$.
\item A semigroup homomorphism between inverse semigroups automatically preserves inverses.
\end{enumerate}

Let $X$ be a topological space, and let $\ipi(X)$ be the set of all the partial homeomorphisms on $X$, by which we mean the homeomorphisms $h:U\to V$ with $U$ and $V$ open sets of $X$. This set has the structure of an inverse semigroup; 
its multiplication is given by composition of partial homeomorphisms wherever this composition is defined: if $h:U\to V$ and $h':U'\to V'$ are partial homeomorphisms then their product is the partial homeomorphism
\[hh':h^{-1}(V\cap U')\to h'(V\cap U')\]
defined at each point of its domain by $(hh')(x)=h'(h(x))$, and the semigroup inversion is the usual inverse of a homeomorphism:
\[h:U\to V\ \ \ \mapsto\ \ \ h^*=h^{-1}:V\to U\;.\]

\begin{definition}\label{def:pseudogroup}
Let $X$ be a topological space. By a \emph{pseudogroup over $X$} will be meant any subset $P\subset \ipi(X)$ which is closed under the multiplication and the inverse of $\ipi(X)$.
A pseudogroup $P$ over $X$ is \emph{full} if it is also closed under identities in the sense that $\ident_U\in P$ for every open set $U\subset X$; and it is \emph{complete} if it is full and for all $h\in\ipi(X)$ and every open cover $(U_\alpha)$ of $\dom(h)$ we have $h\in P$ if $h\vert_{U_\alpha}\in P$ for all $\alpha$.
\end{definition}

The Wagner--Preston theorem (see \cite{Lawson}) asserts that every inverse semigroup is isomorphic to a pseudogroup. However, we shall need more precise terminology:

\begin{definition}\label{def:inversesemigrouprepresentation}
Let $S$ be an inverse semigroup.
By a \emph{representation} of $S$ on a topological space $X$ will be meant a semigroup homomorphism $\rho:S\to\ipi(X)$. The representation is \emph{full} if $\rho$ restricts to an isomorphism $E(S)\to E(\ipi (X))\cong\topology(X)$. By an \emph{inverse semigroup over $X$} will be meant a pair $(S,\rho)$ consisting of an inverse semigroup $S$ equipped with a representation $\rho:S\to\ipi (X)$. If $\rho$ is full (in this case $S$ is necessarily a monoid) then $(S,\rho)$ is said to be a \emph{full inverse semigroup over $X$}. Finally, an inverse monoid representation is called \emph{unital} if it preserves the unit.
\end{definition}

\subsection*{\'Etale groupoids}\label{subsection:EtaleGroupoids}

To start, let us briefly fix notation for groupoids. 
\begin{definition}\label{def:groupoid}
A \emph{(topological) groupoid} $G$ is a pair of topological spaces, the space of \emph{arrows}   $G_1$ and the space of \emph{objects} $G_0$, equipped with continuous maps
\[\xymatrix{
G_2\ar[r]^-m&G_1\ar@(ru,lu)[]_i\ar@<1.2ex>[rr]^r\ar@<-1.2ex>[rr]_d&&G_0\ar[ll]|u
},\]
where $G_2$ is the set $G_1\times_{G_0}G_1$ of \emph{composable pairs} of arrows,
\[G_2=\{(x,y)\in G_1\times G_1\st r(x)=d(y)\}\;,\]
equipped with the subspace topology relative to the product topology on $G_1\times G_1$ (\ie, the pullback of $d$ and $r$ in the category of topological spaces).
The maps $m$, $d$, $r$, $i$, and $u$ are the \emph{multiplication}, \emph{domain},  \emph{range}, \emph{inverse}, and \emph{unit} maps, respectively.
To shorten notation, we shall contract $m(x,y)$ to $xy$, $i(x)$ to $x^{-1}$, and $u(x)$ to $1_x$. We require these maps to satisfy the usual groupoid axioms.
\end{definition}

We shall be mostly concerned with \emph{\'etale} groupoids, in other words with those for which the domain map $d$ is a local homeomorphism, or, equivalently, for which $d$ is open and $u(G_0)$ is open in $G$ (see \cite{Resende} for the latter characterization).

By a \emph{local bisection} of an \'etale groupoid $G$ is usually meant a local section $s:U\to G_1$ of the domain map $d$ on an open set $U\subset G_0$ such that $r\circ s$ is an open embedding of $U$ into $G_0$. Often we shall, usually without any comment, identify the local bisections $s$ with their images $s(U)$, which are the open subsets $V$ of $G_1$ such that the restrictions $d\vert_V$ and $r\vert_V$ are both injective (these sets are called \emph{G-sets} in \cite{Paterson}, following terminology introduced in \cite{Renault}, but we shall avoid this because it clashes with the established terminology for sets equipped with an action of a group $G$).

\begin{definition}
Let $G$ be an \'{e}tale groupoid. The \emph{inverse semigroup of $G$}, $\ipi(G)$, is the set of local bisections of $G$, with multiplication given by pointwise multiplication of (images of) local bisections, and the inverse being similarly calculated pointwise.
\end{definition}

We remark that $\ipi(G)$ acts partially on $G_0$ and in fact is a full inverse semigroup over $G_0$, with the representation $\rho_G:\ipi(G)\to\ipi(G_0)$ defined by
\begin{eqnarray*}
\rho_G(V):d(V)&\stackrel\cong\to& r(V)\\
\rho_G(V)(d(x)) &=& r(x)
\end{eqnarray*}
for each $V\in\ipi(G)$ and $x\in V$.

\begin{definition}\label{definition:GermsFullRepresentation}
Let $(S,\rho)$ be a full inverse semigroup over a topological space $X$. We define the \emph{germ} of $s\in S$ at $x\in\dom(\rho(s))$ to be 
\[\germ_xs=\{t\in S | \exists f\in E(S): ft=fs, \;\; x\in \dom(\rho(t))\cap\dom(\rho(f)) \},\]
and the set of germs of $(S,\rho)$ to be
\[\Germs(S,\rho)=\{(x,\germ_xs) | x\in X, s\in S \}.\]
\end{definition}
We equip $\Germs(S,\rho)$ with the sheaf topology, whose basis consists of open sets of the form, for each $s\in S$,
\begin{eq}\label{Us}
U_s=\{(x,\germ_x s)\st x\in\dom(\rho(s))\}\;.
\end{eq}%
These two constructions can be regarded as inverse to each other, with the connection given by the following theorems (in Section~\ref{section:ACPs} we shall describe conditions on the inverse semigroups under which the two constructions are really inverse to each other):

\begin{theorem}\label{theorem:GermGroupoidFullSemigroup}
Let $(S,\rho)$ be a full inverse semigroup over a topological space $X$. Then the space $\Germs(S,\rho)$ can be given the structure of an \'etale groupoid with object space $X$.
\end{theorem}

\begin{proof}
This construction is an essentially straightforward adaptation of the construction of a local homeomorphism of a sheaf and is done in \cite{Paterson} in a slightly more general setting, namely when $\rho(E(S))$ is only a basis of $X$ rather than the whole topology.

Although also subject to restrictions pertaining to the space $X$, which in \cite{Paterson} has to be Hausdorff, second countable, and locally compact, these restrictions are irrelevant, so we shall describe the construction for full inverse semigroups here.

The groupoid $G=\Germs(S,\rho)$ has the set of arrows
\[G_1=\left\{\left.(x,\germ_xs)\right|s\in S \mbox{ and } x\in\dom(\rho(s))\right\}.\footnote{The arrows are pairs $(x,\germ_x s)$ rather than just the actual germs $\germ_x s$ because unless $X$ is a $T_0$-space we may have $\germ_x s=\germ_y s$ with $x\neq y$ --- see also Section 5.}\]

Note that the projection $d:\Germs(S,\rho)\to X$ defined by $(x,\germ_x s)\mapsto x$ is a local homeomorphism. Also, the subspace topology on $X\subset\Germs(S,\rho)$ coincides with the original topology on $X$. Keeping all these things in mind, we can introduce the topological groupoid structure on the space of germs. Namely, the operations are defined as follows:
 \begin{eqnarray*}
d(x,\germ_x s)&=& x\\
r(x,\germ_x s)&=& \rho(s)(x)\\
1_x&=&(x,\germ_x e)\\
(x,\germ_x s)(\rho(s)(x),\germ_{\rho(s)(x)},t)&=&(x,\germ_x (st))\\
(x,\germ_x s)^{-1}&=&(\rho(s)(x),\germ_{\rho(s)(x)} (s^*))\;.
\end{eqnarray*}
This groupoid is \'etale, for the domain map $d$ is a local homeomorphism as mentioned above. We do not give a direct proof of the soundness of this construction here, since it is similar to that of \cite{Paterson}.
\qed
\end{proof}

Now we prove that every \'etale groupoid arises in this way: 

\begin{theorem}
Let $G$ be an \'etale groupoid, and let $\rho_G:\ipi(G)\to\ipi(G_0)$ be its full representation. Then $\Germs(\ipi(G),\rho_G)\cong G$.
\end{theorem}

\begin{proof}
The standard identification of the total space $E$ of a local homeomorphism $p:E\to X$ with the space of stalks (germs of continuous local sections) of the sheaf of local sections of $p$ gives us a homeomorphism between $\Germs(\ipi(G),\rho_G)$ and $G$, since it is clear that every local section is locally a local bisection: if $x\in U$ and $s:U\to G_1$ is a local section of $d$ then there is an open set $V\subset U$ such that $x\in V$ and $s\vert_V$ is a local bisection. It is then routine to check that the groupoid operations in $\Germs(\ipi(G),\rho_G)$ correspond to those in $G$.
\qed
\end{proof}

\subsection*{Unital representations of inverse monoids}\label{subsection:UnitalRepresentations}

We have seen that any \'etale groupoid $G$ is determined by a full representation of $\ipi(G)$ on the unit space $G_0$, but we shall need to extend this relationship to encompass more general inverse semigroup representations, in particular, in this section we shall see explicitly how any unital representation of an inverse monoid can be turned into a full representation of a larger inverse monoid.

Let $M$ be an inverse monoid, $X$ a topological space, and $\rho:M\to\ipi(X)$ a monoid homomorphism, that is, a unital representation in our terminology.

Define
\[(\Omega(X)\downarrow M)=\{(U,s)\st U\in\Omega(X),\ s\in M,\ U\subset\dom(\rho(s))\}\;,\]
where the pair $(U,s)$ should be thought of as a formal restriction of $s$ to the subspace $U$ of its domain.

\begin{lemma}
The set $(\Omega(X)\downarrow M)$ has a structure of inverse monoid.
\end{lemma}

\begin{proof}
The multiplication can be defined by
\[(U,s)(V,t)=(U\cap\rho(s)^{-1}(V\cap\rho(s)(U)),st)\;,\]
the inverse by
\[(U,s)^*=(\rho(s)(U),s^*)\;,\]
and the unit is \((X,e)\).
The rest is straightforward. \qed
\end{proof}

The monoid $(\Omega(X)\downarrow M)$ is too large for the purposes we have in mind because we shall need the submonoid of idempotents to be isomorphic to $\Omega(X)$, whereas in  $E(\Omega(X)\downarrow M)$ there are in general many copies of each open set, namely $(U,f)$ for each idempotent $f$ of $M$ such that $U\subset\dom(\rho(f))$. Hence, we shall define a quotient of $(\Omega(X)\downarrow M)$ in order to get a full inverse monoid over $X$.

To do this, define an equivalence relation on $(\Omega(X)\downarrow M)$ by
\[(U,s)\sim(V,t)\]
if $U=V$ and there is $f\in E(M)$ such that
$U\subset\dom\rho(f)$ and $fs=ft$.

It is easy to see that this is indeed a congruence relation. We shall denote the congruence class of $(U,s)$ by $[U,s]$ and the quotient $(\Omega(X)\downarrow S)/{\sim}$ by $M_X$.

\begin{lemma}
$M_X$ is an inverse monoid and it is equipped with a full representation $\rho_X: M_X\to\ipi(X)$ given by
\begin{eq}\label{MXrep}
\rho_X([U,s])=\rho(s)\left|_U\right.,\qquad s\in M,\qquad U\subset\dom(\rho(s))\;.
\end{eq}
\end{lemma}
\begin{proof}
If $(U,s)\sim(U,s')$ and $(V,t)\sim(V,t')$ then we have $f, g\in E(M)$ with $fs=fs'$ and $gt=gt'$, and $U$ and $V$ are inside $\dom(\rho(f))$ and $\dom(\rho(g))$ respectively. 

We first claim that the compositions $(U,s)(V,t)$ and $(U,s')(V,t)$ represent the same class in $M_X$. Indeed, $\rho(fs)=\rho(fs')$ implies
\[\rho(s)\left|_U=\rho(s')\left|_U\right.\right.,\]
so that
\[U\cap\rho(s)^{-1}(V\cap\rho(s)(U))=U\cap\rho(s')^{-1}(V\cap\rho(s')(U)),\]
which means that the ``domain'' of the composition is well defined. 

Next, $fst=fs't$ and clearly $f$ is an idempotent for which the domain $U\cap\rho(s)^{-1}(V\cap\rho(s)(U))$ of the compositions under consideration lies inside $U$, which, in turn, is inside $\dom(\rho(f))$ by assumption. 

By a similar computation, now involving $g$ and the appropriate domains, we obtain $(U,s')(V,t)\sim(U,s')(V,t')$ and, together with our previous observation, this leads to $(U,s)(V,t)\sim(U,s')(V,t')$ as desired.

For inverses, take $(U,s)\sim(U,s')$ with $f\in E(M)$ as above. A trivial computation reveals that
\[(s'^*fs')(s^*fs)s^*=(s'^*fs')(s^*fs)s'^*,\]
so that $(U,s)^*\sim(U,s')^*$ using an idempotent $s'^*fs's^*fs$ whose ``domain'' clearly contains $\rho(s)(U)$.

The congruence class $[X,e]$ provides the unit, and so we indeed have an inverse monoid.

Now we prove that the map $\rho_X$ of (\ref{MXrep}) is a representation. It is well defined because $\rho_X[U,s]$ gives us a partial homeomorphism of $X$
by definition. Moreover, it preserves multiplication, inverses and the
unit. Say,
\[
\rho_X( [U,s][V,t] ) = \rho_X[ \rho(s)^{-1}(V \cap \rho(s)(U)), st ] =
\rho(st)|_{\rho(s)^{-1}(V \cap \rho(s)(U))}\;,\]
while
\[\rho_X[U,s] \rho_X[V,t] = \rho(s)|_U \rho(t)|_V\;,\]
and these two partial homeomorphisms of $X$ are the same in $\ipi(X)$.

Further, the representation $\rho_X$ is full, since for every $U\in\Omega(X)$ the class $[U,e]$ gives an idempotent in $M_X$ whose image is $\id_U$.
\qed
\end{proof}

Note that we have the following commuting square of inverse monoid homomorphisms with the vertical homomorphism $\Omega(X)\to M_X$  mapping an open set $U$ to $[U,e]$ and the horizontal homomorphism $M\to M_X$ mapping $s$ to $[\dom(\rho(s)),s]$:
\begin{equation*}
\xymatrix{ E(M)\ar[d]_{\subset}\ar[r]^{\dom\rho} & \Omega(X)\ar[d]\\
M\ar[r] & M_X}
\end{equation*} 

We remark that $M_X$ in general is not a pushout in the category of inverse monoids. However, it has the following universal property.
\begin{lemma}
Let $M'$ be an inverse monoid for which there exist inverse monoid homomorphisms $\alpha$ and $\beta$ making the outer square of the following diagram
\begin{eq}
\vcenter{\xymatrix{ E(M)\ar[d]_{\subset}\ar[r]^{\dom\rho} & \Omega(X)\ar[d]\ar@/^/[ddr]^{\beta} & \\
M\ar[r]\ar@/_/[drr]_{\alpha} & M_X\ar@{.>}[dr]|{\rho'} & \\
& & M'}}\label{equation:Detailed}
\end{eq}%
commutative. Suppose also that 
\[\beta(U\cap\rho(s^{-1})(V\cap\rho(s)(U)))\alpha(s)=\beta(U)\alpha(s)\beta(V)\]
and
\[\alpha(s^*)\beta(U)=\beta(\rho(s)(U))\alpha(s^*)\]
for all $U, V\in\Omega(X)$ and $s\in M$ with $U\subset\dom(\rho(s))$.
Then there exists a unique inverse monoid homomorphism $\rho':M_X\to M'$ as depicted in~\eqref{equation:Detailed}.
\end{lemma}
\begin{proof}
First note that in $M$ we have a decomposition
\[(U,s)=(U,e)(\dom(\rho(s)),s),\qquad U\in\Omega(X), s\in M, U\subset\dom(\rho(s)),\]
which in turn provides us with a similar decomposition for $M_X$:
\begin{equation*}
[U,s]=[U,e][\dom(\rho(s)),s].
\end{equation*}
Thanks to this decomposition, $\rho_X[U,s]$, if it exists, has to be
\[\rho'[U,e]\rho'[\dom(\rho(s)),s]\;,\]
but both of the terms $[U,e]$ and $[\dom(\rho(s)),s]$ come from the preimages $U$ in $\Omega(X)$ and $s$ in $M$ respectively, so that we can define the results of applying $\rho'$ to them using $\beta$ and $\alpha$:
\[\rho'[U,s]=\beta(U)\alpha(s), \qquad U\in\Omega(X), s\in M, U\subset\dom(\rho(s)).\]

Thus we obtain uniqueness of $\rho'$. To check existence, that is, to see that the formula above gives us a homomorphism, we check that
\begin{multline*}
\rho'([U,s][V,t])=\rho'([U\cap\rho(s^{-1})(V\cap\rho(s)(U)),st])=\\
\beta(U\cap\rho(s^{-1})(V\cap\rho(s)(U))\alpha(st)\\
=\beta(U)\alpha(s)\beta(V)\alpha(t)=\rho'([U,s])\rho'([V,t])
\end{multline*}
and 
\begin{multline*}
\rho'([U,s])^*=(\beta(U)\alpha(s))^*=\alpha(s^*)\beta(U)=\\
\beta(\rho(s)(U))\alpha(s^*)=\rho'([\rho(s)(U),s^*])=\rho'([U,s]^*).
\end{multline*}

The commutativity of the resulting diagram comes automatically: starting with $U\in\Omega(X)$, it is being mapped to $[U,e]$ in $M_X$ and then to $\beta(U)\alpha(e)$ in $M'$. But the latter product coincides with $\beta(U)$, for $\alpha(e)$ is the unit of $M'$. For another side of the diagram, an element $s$ in $M$ is being first mapped to $[\dom(\rho(s)),s]$ in $M_X$ and after that via $\rho'$ to $\beta(\dom(\rho(s)))\alpha(s)$. Consider an idempotent $ss^*$ whose ``domain'' clearly covers the one of $s$ (in fact, they are the same) and trace it through the outer square in~\eqref{equation:Detailed} to obtain
\[\beta(\dom(\rho(s)))=\beta(\dom(\rho(ss^*)))=\alpha(ss^*).\]
Multiplying this identity by $\alpha(s)$ on the right, we get $\beta(\dom(\rho(s)))\alpha(s)=\alpha(s)$ as required.\qed
\end{proof}

Of course, a pushout in the category of inverse monoids also exists, but any concrete description of such a pushout would require adding new idempotents $sfs^*$ for $s\in M$ and $f$ an ``old'' idempotent in $M$.

\subsection*{Germ groupoids revisited}

The germ groupoid of $(M_X,\rho_X)$ from the previous section has a direct description in terms of the germ groupoid of $(M,\rho)$, as we shall now see. Further, we exhibit this correspondence in an even more general context.

\begin{definition}
Let $(M,\rho)$ be an inverse monoid over a topological space $X$ (as usual, we assume that $\rho$ is a unital representation). As in Definition~\ref{definition:GermsFullRepresentation}, define the \emph{germ} of $s\in M$ at $x\in\dom(\rho(s))$ to be 
\[\germ_xs=\{t\in M | \exists f\in E(M): ft=fs, \; x\in \dom(\rho(t))\cap\dom(\rho(f)) \}\]
and the space of germs to be
\[\Germs(M,\rho)=\{(x,\germ_xs) | x\in X, s\in M \}.\]
Further, we topologize the space of germs $\Germs(M,\rho)$ by making a topology basis out of the sets
\[V_{s,U}=\left\{ (x, \germ_xs) | x\in U \right\}\]
for $s\in M$ and $U\subset\dom(\rho(s))$, $U\in\Omega(X)$. Note that this sheaf topology is different from the one introduced in Theorem~\ref{theorem:GermGroupoidFullSemigroup} before.
\end{definition}

We have a straightforward counterpart of Theorem~\ref{theorem:GermGroupoidFullSemigroup}:

\begin{theorem}\label{theorem:GermGroupoidUnitalMonoid}
Let $(M,\rho)$ be a unital inverse monoid representation over a topological space $X$. Then the space $\Germs(M,\rho)$ can be given the structure of an \'etale groupoid with object space $X$.
\end{theorem}

We shall not give a detailed proof of this theorem here, but mention that one defines all the structure maps exactly as in Theorem~\ref{theorem:GermGroupoidFullSemigroup} before and thus obtains a groupoid with object space $X$. The only nontrivial part is to show that this groupoid is \'etale. We shall see that the groupoid in question is the same as the one for a full representation $\rho_X$ obtained from $\rho$ and thus reduce it to the case already discussed.

\begin{theorem}\label{theorem:GermGroupoidsEquivalenceMonoids}
Let $(M,\rho)$ be a unital inverse monoid representation over a topological space $X$ and $(M_X,\rho_X)$ --- an induced full inverse monoid as constructed in the previous section. Then \[\Germs(M,\rho)\cong\Germs(M_X,\rho_X).\]
\end{theorem}
\begin{proof}
We start by writing explicitly what germs in $\Germs(M_X,\rho_X)$ are. Take $[U,s]\in M_X$ and $x\in\dom(\rho_X([U,s])=\dom(\rho(s)|_U)=U$. Then
\begin{multline}\label{equation:GermsMX}
\germ_x[U,s]=\left\{[V,t]\in M_X \,|\, \exists [F,f]\in E(M_X)\right.\\
\left.\mbox{ with } x\in F,\; [F\cap U,fs]=[F\cap V,ft] \right\}.
\end{multline}
In order to identify this germ with $\germ_xs\in\Germs(M,\rho)$, we check that $t$ from $[V,t]$ above also belongs to $\germ_xs$. This is true, for clearly $x\in V\subset\dom(\rho(t))$ and $gfs=gft$ for some $g\in E(M)$ covering $F\cap U\ni x$. This shows that the mapping 
\[\Germs(M_X,\rho_X)\to\Germs(M,\rho): \quad \germ_x[U,s]\mapsto\germ_xs\]
is surjective and well defined. To check the injectivity of such a correspondence, take $\germ_x[U,s]$ and $\germ_x[V,s]$ which are both mapped to $\germ_xs$. In fact, for an idempotent $[U\cap V,e]$ of $M_X$, its domain under $\rho_X$ contains $x$ and one has
\[[U\cap V,e][U,s]=[U\cap V,e][V,s],\]
which means that the germs are indeed the same, according to~\eqref{equation:GermsMX}.

Next, we clarify why the groupoid structure maps in $\Germs(M,\rho)$ and $\Germs(M_X,\rho_X)$ are the same. This is so because in the correspondence between germs established above the domains are $x$ and the ranges are $\rho(s)(x)=\rho_X([U,s])(x)$. The inverses are $\germ_{\rho(s)(x)}s^*$ and $\germ_{\rho_X([U,s])(x)}[\rho(s)U,s^*]$ which clearly correspond to each other, the units $1_x$ are defined using the germs of $e$ and $[X,e]$ respectively, which also correspond to each other. Finally, the composition is preserved as well: for two composable arrows $\germ_x[U,s]$ and $\germ_{\rho(s)(x)}[V,t]$ in $\Germs(M_X,\rho_X)$ their composition is
\[\germ_x[U\cap\rho(s)^{-1}(V\cap\rho(s)(U)),st],\]
which corresponds to $\germ_x(st)$ in $\Germs(M,\rho)$ as expected.

We finish by checking that the topologies imposed on $\Germs(M,\rho)$ and $\Germs(M_X,\rho_X)$ are compatible. An elementary open set in the latter space, 
\[V_{[U,s]}=\{ \germ_x[U,s]\,|\,x\in\dom(\rho_X([U,s]))=U\,\},\] 
corresponds to $V_{s,U}$ in the former one, and vice-versa, which completes the proof.
\qed
\end{proof}

\subsection*{Wide representations of inverse semigroups}

Finally, we want to treat the following slightly more general case.

\begin{definition}
A representation $\rho: S\to\ipi(X)$ of an inverse semigroup $S$ on a topological space $X$ will be called \emph{wide} if $\rho(E(S))$ covers $E(\ipi(X))$. We shall refer to such $S$ as  a \emph{wide inverse semigroup over $X$}.
\end{definition}

In particular, any full representation is wide. More generally, any inverse monoid representation is a wide representation. For a rather extreme example, the  Wagner--Preston representation of any inverse semigroup is wide.

Given a wide inverse semigroup $(S,\rho)$ over a topological space $X$, define the germ of $s\in S$ at $x\in\dom(\rho(s))$ and the space of germs exactly as above to be 
\[\germ_xs=\{t\in S | \exists f\in E(S): ft=fs, \; x\in \dom(\rho(t))\cap\dom(\rho(f)) \}\]
and
\[\Germs(S,\rho)=\{(x,\germ_xs) | x\in X, s\in S \}\]
respectively. The topology on the space of germs is generated by the basis consisting of
\[V_{s,U}=\left\{ (x, \germ_xs) | x\in U \right\}\]
for $s\in S$ and $U\subset\dom(\rho(s))$, $U\in\Omega(X)$, as for the unital monoid representation above.

One can extend Theorem~\ref{theorem:GermGroupoidUnitalMonoid} to

\begin{theorem}\label{theorem:GermGroupoidWideSemigroup}
For a wide inverse semigroup $(S,\rho)$, the space $\Germs(S,\rho)$ can be given a structure of \'etale groupoid with object space $X$.
\end{theorem}

Again, our goal is not to prove this theorem independently, but rather by means of connecting it with other constructions we have studied.

\begin{theorem}\label{theorem:GermGroupoidsEquivalenceWide}
Let $(S,\rho)$ be a wide inverse semigroup over a topological space $X$ and the inverse monoid $S_e$ be the result of adjoining a unit to $S$. Further, extend $\rho$ to a monoid representation $\rho_e$:
\[\xymatrix{ S \ar[r]^{\subset}\ar[dr]_{\rho} & S_e \ar@{.>}[d]^{\rho_e} \\ & \ipi(X)}\]
Then \[\Germs(S,\rho)\cong\Germs(S_e,\rho_e).\]
\end{theorem}
\begin{proof}
Since we impose the same topology on the space of germs in both cases, all we need to check is that the spaces of germs are the same as sets. And this is indeed true, because $\rho$ is wide and thus we just add the unit $e$ to the germs that contain idempotents (this will ``enrich'' some germs, but it will not add any new ones). \qed
\end{proof}

Note that for $\Germs(S_e,\rho_e)$ all the units can be written as germs of $e$ as in the full representation case which was not the case right away for wide representations.

\begin{corollary}\label{corollary:GermGroupoidWide}
Any wide inverse semigroup $(S,\rho)$ over a topological space $X$ determines an inverse monoid with a full representation on $X$, and their germ groupoids are isomorphic.
\end{corollary}

We finish this section by mentioning that the notion of a wide inverse semigroup is more general than the localizations of \cite{Kumjian,Paterson}, since we do not impose any conditions on the topology of $X$ and also do not require the idempotents to provide a basis for the topology, but just a cover.

\section{Example: the universal groupoid of Paterson}\label{section:UniversalGroupoid}

We shall show that the \emph{universal groupoid of an inverse semigroup} introduced in~\cite{Paterson} can be obtained using our wide representation techniques. This section is intended to be as self-contained as possible, so we start by describing the construction and then comment on the properties of the universal groupoid and some motivations behind it.

Let $S$ a countable inverse semigroup. Consider the set
\[X=\left\{\mbox{non-zero multiplicative functions }x: E(S)\to\{0, 1\}\right\}.\]
For each $s$ in $S$ let $D_s$ denote the subset of $X$ consisting of all the functions $x$ for which $x(ss^*)=1$. We topologize $X$ by specifying as a basis for the topology the collection of all the sets
\[D_f\cap(X\backslash D_{f_1})\cap(X\backslash D_{f_2})\cap\dots\cap(X\backslash D_{f_n})\]
where $f, f_1, f_2, \dots, f_n$ are idempotents of $S$ such that $ff_i=f_i$ for all $i=1,\ldots n$.
In this way $X$ becomes a locally compact totally disconnected space (and hence Hausdorff).

Now we construct a representation $\rho_u$ of $S$ on $X$. For $s\in S$, we define  a homeomorphism
\[\rho_u(s): D_s \to D_{s^*}\]
by sending $x\in D_s$ to $y\in D_{s^*}$ defined by $y(f)=x(sfs^*)$.

It is easy to see that this representation is wide: for any $x\in X$ there ought to be an idempotent $f$ with $x(f)=1$, since $x$ is not identically $0$. Then $\rho_u(f)$ is a partial identity of $X$ with domain $D_f$ which tautologically contains $x$. Therefore we can construct the germ groupoid $G_u$ corresponding to $S$ and $\rho_u$ as prescribed in the previous section. This groupoid is called the \emph{universal groupoid} of $S$. It enjoys the following properties:

\begin{enumerate}
\item $C^*(S)=C^*(G_u)$.
\item $C^*_{red}(S)=C^*_{red}(G_u)$.
\item For any other ample $S$-groupoid $G$ (see \cite{Paterson}), its unit space $G_0$ is homeomorphic to a closed invariant subspace $Y$ of $(G_u)_0$, and there exists a continuous open surjective $S$-equivariant homomorphism
\[\phi: G_u|_Y\to G, \;\; \phi|_Y=id.\]
\end{enumerate}

We comment that our construction, unlike the original one of Paterson in~\cite{Paterson}, does not appeal to the localisation techniques, and therefore does not require the inverse semigroup $S$ to be countable (and the germ groupoid construction from the previous section does not require the space $X$ to be locally compact Hausdorff in general). Of course, for the proof of the analytic properties of the universal groupoid listed above, these constrains are still required.

\section{Example: the translation groupoid of Skandalis, Tu, and Yu}\label{section:TranslationGroupoid}
\subsection*{The Stone--\v{C}ech compactification via ultrafilters}

We start by briefly reviewing some elements of the construction of the Stone--\v{C}ech compactification of a discrete space, mostly to fix the notation. For more details the reader is referred to~\cite{Johnstone}.

\begin{definition}
Let $X$ be a space with the discrete topology. A collection $\mathcal F$ of subsets of $X$ is called a \emph{filter} if the following conditions are satisfied:
\begin{itemize}
\item $\emptyset\notin\mathcal F$.
\item Whenever $U\subset V$ for some subsets $U, V$ of $X$, and $U\in\mathcal F$, then  $V\in\mathcal F$.
\item For $U, V\in\mathcal F$, $U\cap V\in\mathcal F$.
\end{itemize}
\end{definition}

One can identify the Stone--\v{C}ech compactification of a discrete space $X$ with the space of \emph{ultrafilters}, which are by definition maximal (with respect to inclusion) filters of subsets of $X$. The topology on this space is generated by the sets
\[\tilde U=\{ \mathcal F | U\in\mathcal F\}, \qquad U\subset X.\]

The embedding of $X$ into $\beta X$ for our discrete space case in this model is done by  mapping $x\in X$ to $\mathcal F_x$, the \emph{principal} filter at $x$ which is defined to consist of all subsets of $X$ which contain $x$. 

Further, if one identifies $X$ with the image of its embedding into $\beta X$, any basic open set $\tilde U$ of $\beta X$ has $U$ as its trace on $X$. Moreover, given a set $U$ in $X$, its closure in $\beta X$ is precisely $\tilde U$.

We conclude this review with the following technical result. 
\begin{lemma}\label{lemma:HomeomorphismBetaU}
Let $U$ be a subset of a discrete space $X$. Then there exists a canonical homeomorphism between $\beta U$ and $\tilde U$ (the latter one being a subset of the ambient space $\beta X$ with the induced topology.)
\end{lemma}
\begin{proof}
If we take some ultrafilter $\mathcal F$ in $\tilde U$ then the family
\[\mathcal F_U:=\{ F\cap U | F\in\mathcal F \}\]
is an ultrafilter in $U$. Conversely, given an ultrafilter $\mathcal G$ of subsets of $U$, one can formally define 
\[\mathcal G_X=\{ G\cup A | G\in\mathcal G, \; A\subset X\backslash U\},\]
which turns out to be an ultrafilter on $X$ extending $U$. This correspondence between ultrafilters on $U$ and the ones on $X$ within $\tilde U$ is bijective (the `shrinking' and `enlarging' procedures are inverse to each other).

Under this correspondence the standard basis of the topology on $\beta U$, namely the one comprised out of the sets
\[\tilde A=\{\mathcal G\in\beta U | A\in\mathcal G\}, \quad A\subset U,\]
corresponds (element-wise) to
\begin{eq}\label{eqn:InducedBasis}\{\mathcal F\in\beta X | A\in\mathcal F\}.\end{eq}%
We claim that the latter sets are the intersections of the elements of the standard basis for $\beta X$ with $\tilde U$. Indeed, given $V\subset X$, $\tilde U\cap\tilde V=\widetilde{U\cap V}$, so that taking $A=U\cap V$ the condition that $V$ runs over all subsets of $X$ is equivalent to that of $A$ running over all subsets of $U$. This shows that the sets in~\eqref{eqn:InducedBasis} form a basis for $\tilde U$ in $\beta X$ and, moreover, the correspondence which we established respects the aforementioned bases of the topologies, hence it is a homeomorphism.
\qed\end{proof}

\subsection*{The Stone--\v Cech compactification of the unit space}

Having a partial homeomorphism $h: U\to V$ between two subsets $U, V$ of $X$, we can extend it to a homeomorphism $\beta h: \beta U\to\beta V$ or, equivalently, by the virtue of Lemma~\ref{lemma:HomeomorphismBetaU}, to a homeomorphism $\tilde h: \tilde U\to\tilde V$.

Now we are in position to prove the Stone--\v Cech extension theorem for discrete groupoids.

\begin{theorem}\label{theorem:StoneCechExtensionDiscrete}
Let $G$ be a discrete groupoid. Then there exists an \'etale groupoid (which we shall denote by $\beta_0G$) such that
$G$ is a subgroupoid of $\beta_0G$ and $(\beta_0G)_0=\beta G_0$.
\end{theorem}
\begin{proof}
The inverse semigroup $\ipi(G)$ has a representation $\rho_G$ on $G_0$ (notice that all the partial homeomorphisms in this presentation are just partial bijections.) We extend it to $\widetilde{\rho_G}: \ipi(G)\to\ipi(\beta G_0)$ by means of taking every partial homeomorphism $h: U\to V$ on $G_0$ from $\rho_G(\ipi(G))$ and extending it to a partial homeomorphism $\tilde h: \tilde U\to \tilde V$ on $\beta G_0$. It is easy to see that by doing such an extension we indeed obtain a representation: $(\tilde h)^{-1}=\widetilde{(h^{-1})}$ and $\tilde h\tilde h'=\widetilde{hh'}$.

Notice that whilst the representation $\rho_G$ is full (so that $\rho_G(E(\ipi(G)))\cong\Omega(G_0)$), $\widetilde{\rho_G}$ is wide, due to the fact that  $\widetilde{\rho_G}(E(\ipi(G)))$ contains all the sets $\tilde U$ with $U\in\Omega(G_0)$, which form a basis for $\Omega(\beta G_0)$ and hence cover it.

The statement now follows from a direct application of Corollary~\ref{corollary:GermGroupoidWide}. To see that $G$ is a subgroupoid of $\beta_0G$, recall that $G$ can be identified with the germ groupoid $\Germs(\ipi(G),\rho_G)$. For any $x\in G_0$ and any $s\in\ipi(G)$ with $x\in\dom(\rho_G(s))$, the germ of $s$ at $x$ as an   arrow in $G$ can be viewed as an arrow in $\beta_0G$, since the representation $\widetilde{\rho_G}$ coincides with $\rho_G$ on $X$ (more specifically, for $x\in X$, the conditions that $x$ belongs to the domain of $\widetilde{\rho_G}(t)$ and to the domain of $\rho(t)$ for some $t\in\ipi(G)$, are equivalent. But aside from domain conditions, the definitions of the germs forming $G$ and $\beta_0G$ are the same.)
\qed\end{proof}

\begin{remark}
In general $(\beta_0G)_1\ne\beta G_1$, so that the newly constructed groupoid is not the Stone--\v Cech compactification of the original groupoid $G$.
\end{remark}

\subsection*{Digression: the translation groupoid}

We shall discuss one more specialized case of the germ groupoid construction and the Stone--\v Cech extension when the inverse semigroup is a pseudogroup over a discrete topological space.
The motivation for such a digression is that Skandalis, Tu, and Yu have proven that the coarse Baum--Connes conjectures for the discrete coarse space and the resulting groupoid (which they called the \emph{translation groupoid}) are equivalent.

Here we present a brief account of the ideas involved. For more details on the construction and the ambient context consult the original paper~\cite{SkandalisTuYu} of Skandalis, Tu, and Yu; a more comprehensive account on coarse geometry can be found in~\cite{Roe}.

Let $X$ be an infinite set endowed with the discrete topology. In what follows, it is convenient to regard $X\times X$ as a \emph{pair groupoid}, that is, the groupoid with object space $X$, arrow space $X\times X$, $d$ and $r$ being the first and second projections $\pi_1,\pi_2:X\times X\to X$, etc.
\begin{definition}
A \emph{coarse structure} on $X$ is a collection $\mathcal E$ of nonempty subsets of $X\times X$, called \emph{controlled sets}, such that every singleton of $X\times X$ belongs to $\mathcal E$ and $\mathcal E$ is closed with respect to taking
\begin{itemize}
\item Subsets;
\item Finite unions;
\item Inverses (forming a new set consisting of the inverses of the elements of the original set in the pair groupoid sense);
\item Products (forming a new set out of the products of all composable elements from two controlled sets);
\end{itemize}
Such $X$ together with a coarse structure is called a \emph{coarse space}.
\end{definition}

One important example of coarse spaces comes from a metric. Starting from a metric space $X$, one can define the coarse structure to contain all sets
\[E\subset X\times X \mbox{ such that } \exists N: \forall(x,y)\in E \; \dist(x,y)<N.\]
Not all coarse structures come from a suitable metric, however; in some important cases one can in principle impose different, yet equivalent metrics on the same space (the standard example is the word metric on a finitely generated group with respect to different choice of the generating set; all metrics in this example are bi-Lipschitz equivalent) --- it turns out that the resulting coarse structures are the same. Thus the coarse space approach allows one to study finitely generated groups as metric spaces without explicit reference to a particular generating set.

When the diagonal $\{(x,x)|x\in X\}$ of a coarse space $(X,\mathcal E)$ is controlled, the coarse structure is called \emph{unital}. This happens, for instance, for coarse structures which arise from a metric.

Given a coarse space $(X,\mathcal E)$, the set $S=\ipi(X)\cap\mathcal E$ of controlled partial bijections on $X$ is a pseudogroup with a naturally defined representation on $X$. The subtle issue is that whilst this representation is indeed wide, it does not have to be full, and therefore its extension to a representation on $\beta X$ by extending all partial bijections to the Stone--\v Cech compactification as in the previous subsection is not necessarily wide.

In the case where the original coarse structure $\mathcal E$ comes from a metric, it is unital and therefore every subset $E$ of the diagonal of $X\times X$ is controlled.  Since for each such subset both coordinate projections are injective, $E$ belongs to $E(S)$ and, being viewed as a partial identity on $X$, it extends to an idempotent on $\beta X$. It is clear that for every subset $U$ of $X$ the identity on $U$ can be represented by such a controlled set $E$, and the sets $\tilde U$ which are the domains of the extended idempotents form a basis for $\beta X$, and so the extended representation is wide. By applying Corollary~\ref{corollary:GermGroupoidWide} we can produce an \'etale groupoid $G(X)$, the \emph{translation groupoid of $X$}, with the unit space being the Stone--\v Cech compactification of the space X. 

Notice that the original representation of $S$ on $X$ is wide (even for the general nonunital case), and this allows us to construct the germ groupoid for it right away. In fact, since each germ in the discrete topology is simply a singleton bijection $\{(x, y)\}$, the resulting germ groupoid is the pair groupoid $X\times X$. For the unital case this means that $X\times X$ is a subgroupoid of $G(X)$, for the natural representation of $S$ is a subrepresentation of the extended one.

\section{Complete inverse semigroups}\label{section:ACPs}

In this section we shall close the circle by providing a characterization of the inverse semigroups of the form $\ipi(G)$ for \'etale groupoids $G$, thereby establishing an equivalence (non functorial) between \'etale groupoids over $X$ and full representations of such inverse semigroups over $X$. In addition we shall provide a brief account of the relation between these results and those of \cite{Resende} concerning localic groupoids and quantales.

\subsection*{Characterization of the monoids of local bisections}

In what follows, we shall make use of the \emph{natural order} on an inverse semigroup. Given an inverse semigroup $S$, this is a partial order and it is defined as follows:
\[s\le t\iff s=ft\textrm{ for some }f\in E(S)\;.\]
Further, the \emph{join} of a subset $Z\subset S$ is an element $\V Z$ of $S$ which is the least upper bound of the elements of $Z$. For more details on these notions and the relevant discussion refer to~\cite{Lawson}. 

\begin{definition}
Let $S$ be an inverse semigroup. Two elements $s,t\in S$ are said to be \emph{compatible} if both
$st^*$ and $s^*t$ are idempotents. A subset $Z\subset S$ is \emph{compatible} if any two elements in $Z$ are compatible. Then $S$ is said to be \emph{complete} if every compatible subset $Z$ has a join $\V Z$ in $S$ (hence, $S$ is necessarily a monoid with $e=\V E(S)$). 
\end{definition}

We are interested in the following class of inverse semigroups.

\begin{definition}
By a \emph{complete inverse semigroup over a space $X$} is meant a complete inverse semigroup $S$ equipped with a full representation $S\to\ipi(X)$.
\end{definition}

The existence of the full representation in this definition has an important consequence for $S$, namely the semilattice of idempotents $E(S)$ is isomorphic to the topology of a space and thus it is a \emph{locale} (see \cite{Johnstone}). It follows (see \cite{Lawson}) that the inverse semigroup is \emph{infinitely distributive} in the sense that for all compatible subsets $Z\subset S$ and all $s\in S$ the set $sZ$ is compatible and we have
\[s \V Z = \V(sZ)\;.\]

Another important fact related to joins, which in particular implies that any full representation of a complete inverse semigroup preserves joins of compatible sets, is that a homomorphism of complete inverse semigroups
$h:S\to T$ whose restriction $h\vert_{E(S)}:E(S)\to E(T)$ preserves arbitrary joins necessarily preserves joins of all the compatible sets \cite[Proposition 2.10-3]{Resende}; that is, for all compatible sets $Z\subset S$ the image set $h(Z)$ is compatible and we have $\V h(Z) = h\left(\V Z\right)$.


We are now ready to give a characterization of the inverse semigroups that arise from \'{e}tale groupoids:

\begin{theorem}
Let $G$ be an \'{e}tale groupoid with unit space $X$. Then $(\ipi(G),\rho_G)$ is a complete inverse semigroup over $X$. Any complete inverse semigroup over $X$ arises in a similar way from an \'{e}tale groupoid with unit space $X$.
\end{theorem}

\begin{proof}
It is easy to see that $\ipi(G)$ is a complete inverse semigroup. For the converse, let $(S,\rho)$ be an arbitrary complete inverse semigroup over $X$, and let $G=\Germs(S,\rho)$. We shall show that $S$ and $\ipi(G)$ are isomorphic, much in the same way in which one shows that a sheaf is isomorphic to the sheaf of local sections of its local homeomorphism. First let us
consider the map $s\mapsto U_s$, defined as in (\ref{Us}):
\[U_s=\{(x,\germ_x s)\st x\in\dom(\rho(s))\}\;.\]
This assignment clearly is a semigroup homomorphism $S\to\ipi(G)$, and it is injective due to infinite distributivity, for if $\dom(\rho(s))=\dom(\rho(t))$ (equivalently, $ss^*=tt^*$) then the condition
$\germ_x s=\germ_x t$ for all $x\in\dom(\rho(s))$ implies that there is a cover $(f_x)$ of $ss^*$ such that for each $x\in\dom(\rho(s))$ we have $f_x s=f_x t$, and thus
\[s=ss^*s=(\V_x f_x)s=\V_x (f_x s)=\V_x (f_x t)=(\V_x f_x) t=tt^*t=t\;.\]

Now let $U$ be an open subset of $G$ such that both the domain and range maps are injective when restricted to $U$. Then, by definition of the topology of $G$, $U$ is a union of sets $U_s$. Let $U_s$ and $U_t$ be two such sets. For all $x\in\dom(\rho(s))\cap\dom(\rho(t))$ we must have a unique arrow of $G$ in $U$ with domain $x$. But $U_s\cup U_t\subset U$, and thus both $(x,\germ_x s)$ and $(x,\germ_x t)$
belong to $U$, therefore implying that $\germ_x s=\germ_x t$; that is, there is an idempotent $f_x\le ss^*tt^*$ such that $x\in\dom(\rho(f_x))$ and $f_x s= f_x t$, and thus
\[(ss^*tt^*) s=(\V_x f_x) s=\V_x (f_x s)=\V_x (f_x t)=(ss^*tt^*) t\;.\]
Hence we have $s^*t\in E(S)$. Similarly, considering any point
$x\in\cod(\rho(s))\cap\cod(\rho(t))$ we conclude, because there must be a unique element in $U$ with codomain $x$, that
$s(s^*st^*t)=t(s^*st^*t)$ [this is immediate from the previous argument because $x\in\dom(\rho(s^*))\cap\dom(\rho(t^*))$], and thus $st^*\in E(S)$. We have thus proved that the set $Z$ that indexes the cover $U=\bigcup_{s\in Z} U_s$ is compatible. Since $S$ is complete, we have a join $\V Z$ in $S$, and it is now clear that $U_{\V Z}=U$, for
\[U_{\V Z}=\bigcup_{s\in Z}U_s=U\;,\]
where we have used the fact that the assignment $U_{(-)}:s\mapsto U_s$ preserves all the joins of compatible sets, which is a consequence of the fact that its restriction to the idempotents does, since that restriction is an isomorphism. Finally, we obviously have $\rho(s)=\rho_G(U_s)$ for all $s\in S$; that is, $U_{(-)}$ commutes with the representations $\rho$ and $\rho_G$, and thus $(S,\rho)$ and $(\ipi(G),\rho_G)$ are the same up to isomorphism. \qed
\end{proof}

Hence, we have arrived at a bijection:

\begin{corollary}\label{firstequiv}
Let $X$ be a topological space. The notions of complete inverse semigroup over $X$ and of \'etale groupoid with unit space $X$ are equivalent up to isomorphisms.
\end{corollary}

\subsection*{Localic germ groupoids}

It is easy to see that if the space $X$ is sober (\ie, $X$ is homeomorphic to the spectrum of a locale --- equivalently, the assignment $x\mapsto\overline{\{x\}}$ from $X$ to the set of irreducible closed subsets of $X$ is a bijection \cite{Johnstone}) then any full representation $\rho:S\to\ipi(X)$ of a complete inverse semigroup is uniquely determined by the isomorphism $E(S)\cong\topology(X)$ because it follows from a representation by conjugation on the locale $E(S)$: each $s\in S$ determines an isomorphism of open sublocales \[\downsegment(ss^*)\to\downsegment(s^*s)\] whose inverse image homomorphism sends each $f\in\downsegment(s^*s)$ to $sfs^*$. Hence, \ref{firstequiv} restricts to an equivalence that generalizes in a nice way, albeit non-functorially, the well known duality between spatial locales and sober spaces (see \cite{Johnstone}). In order to state it, let us follow the terminology of \cite{Resende} and call any complete and infinitely distributive inverse semigroup an \emph{abstract complete pseudogroup (ACP)}, and let us say that an ACP is \emph{spatial} if its locale of idempotents is spatial. We shall abbreviate $\Germs(S,\rho)$ to $\Germs(S)$.

\begin{corollary}\label{secondequiv}
The notions of spatial ACP and of sober \'etale groupoid (an \'etale groupoid whose unit space is sober) are equivalent: if $S$ is an ACP and $G$ is a sober \'etale groupoid then we have isomorphisms
\begin{eqnarray*}
G&\cong&\Germs(\ipi(G))\\
S&\cong&\ipi(\Germs(S))\;.
\end{eqnarray*}
\end{corollary}
This correspondence is a particular instance of the more general correspondence between localic \'etale groupoids and ACPs in \cite{Resende}, which does not depend on spatiality and uses quantales as a mediating structure:
\begin{enumerate}
\item If $S$ is an ACP then its full join completion, denoted by $\lcc(S)$, is a quantale of a kind known as inverse quantal frame;
\item Any inverse quantal frame $Q$ determines an associated localic \'etale groupoid $\groupoid(Q)$;
\item Any localic \'etale groupoid $G$ has an associated ACP $\ipi(G)$;
\item If $S$ is an ACP then $S\cong\ipi(\groupoid(\lcc(S)))$;
\item If $Q$ is an inverse quantal frame then $Q\cong\lcc(\ipi(\groupoid(Q)))$;
\item If $G$ is a localic \'etale groupoid then $G\cong\groupoid(\lcc(\ipi(G)))$.
\end{enumerate}

\begin{theorem}
If $S$ is a spatial ACP then $\groupoid(\lcc(S))$ is spatial and its spectrum is homeomorphic to $\Germs(S)$.
\end{theorem}

\begin{proof}
If $S$ is a spatial ACP it is easy to verify that both $\lcc(S)$ and $\topology(\Germs(S))$ are inverse quantal frames whose associated inverse semigroups of partial units (see \cite{Resende}) are isomorphic to $S$. Hence, the two quantales $\lcc(S)$ and $\topology(\Germs(S))$ are isomorphic and the intended result follows.
\qed
\end{proof}
This explains how the results in \cite{Resende} can be regarded as being a generalization of those of \ref{secondequiv} in the sense of allowing one to define the notion of ``germ groupoid'' in the absence of spatiality, and in toposes beyond the category of sets: $\groupoid(\lcc(S))$ is the required generalization of $\Germs(S)$. It also suggests localic versions of the examples seen in this paper: a suitable generalization of the universal groupoid of an inverse semigroup is likely to involve the patch construction for locales \cite{Escardo}, whereas translation groupoids should rely on Stone-\v{C}ech compactification for locales, for instance as described in \cite{Johnstone}.

We conclude with a simple observation regarding germs and the points of inverse quantal frames:

\begin{corollary}
Let $S$ be an ACP. The germs of $S$ can be identified with the filters $F$ of $S$ that are compatibly prime in the sense that, for all compatible sets $Z\subset S$, if $\V Z\in F$ then $z\in F$ for some $z\in Z$.
\end{corollary}

\begin{proof}
The germs of $S$ correspond to the locale points of $\lcc(S)$, which are the homomorphisms of locales $\lcc(S)\to 2$, where $2$ is the two element chain. The universal property of the principal ideal embedding $S\to\lcc(S)$ \cite{Resende} identifies the points with the non-zero maps $p:S\to 2$ that preserve binary meets and joins of compatible sets, and thereby with the subsets $p^{-1}(1)$, which are the intended filters. \qed
\end{proof}
By direct computation it can be verified that this identification is strict: the germs, which are subsets of $S$, are precisely the compatibly prime filters.


~\\
{\sc
Centro de An\'alise Matem\'atica, Geometria e Sistemas Din\^amicos
Departamento de Matem\'{a}tica, Instituto Superior T\'{e}cnico\\
Universidade T\'{e}cnica de Lisboa\\
Av.\ Rovisco Pais 1, 1049-001 Lisboa, Portugal}\\
{\it E-mail:} {\sf $\{$matsnev,pmr$\}$@math.ist.utl.pt}


\begin{thebibliography}{8}
{

\bibitem{Escardo} M.H.\ Escard\'o, The regular-locally-compact coreflection of a stably locally compact locale, \textit{J.\ Pure Appl.\ Algebra} 157 (2001) 41--55.

\bibitem{Johnstone}
P.T.\ Johnstone,
\textit{Stone Spaces},
Cambridge Stud.\ Adv.\ Math.,
vol.\ 3, Cambridge Univ.\ Press, 1982.

\bibitem{Kumjian}
A.\ Kumjian,
On localizations and simple C*-algebras,
\textit{Pacific J.\ Math.} 112 (1984) 141--192.

\bibitem{Lawson}
M.V.\ Lawson,  \textit{Inverse Semigroups --- The Theory of Partial Symmetries}, World Scientific, 1998.

\bibitem{Paterson}
A.L.T.\ Paterson,
\textit{Groupoids, Inverse Semigroups, and Their Operator Algebras},
Birkh\"{a}user, 1999.

\bibitem{Renault}
J.\ Renault, \textit{A Groupoid Approach to C*-algebras}, Lect.\ Notes Math.\ 793, Springer-Verlag, 1980.

\bibitem{Resende}
P.\ Resende,
\'{E}tale groupoids and their quantales, \textit{Adv.\ Math.} 208 (2007) 147--209.

\bibitem{SkandalisTuYu} G.\ Skandalis, J.-L.\ Tu, G.\ Yu, The coarse Baum--Connes conjecture and groupoids, \textit{Topology} 41 (2002) 807--834.

\bibitem{Roe} J. Roe, \textit{Lectures on Coarse Geometry}, University Lecture Series, American Math.\ Soc., 2003.
}
\end{thebibliography}
\end{document}